\documentclass{amsart}
\newtheorem*{theorem}{\rm THEOREM}
 \newtheorem{lemma}{\rm LEMMA}

 \theoremstyle{remark}

 \numberwithin{equation}{section}

\begin{document}

\title{SOME REMARKS RELATED TO MAEDA'S CONJECTURE}
\author{M. Ram Murty \& K. Srinivas}
\address{Department of Mathematics and Statistics, Queen's University, Jeffery Hall, 99 University Avenue, Kingston, Ontario, Cananda K7L 3N6}
\address{The Institute of Mathematical Sciences, CIT Campus, Taramani, Chennai, 600 113, Tamilnadu, India}
\email[M. Ram Murty]{murty@mast.queensu.ca}
\email[K. Srinivas]{srini@imsc.res.in} 
\begin{abstract}
In this article we deal with the problem of counting the number of pairs of normalized eigenforms $ (f,g) $ of weight $k$ and level $N$ such that $ a_p (f) = a_p (g) $ where $a_p (f) $ denotes the $p-$th Fourier coefficient of $f$. Here $p$ is a fixed prime.  
\end{abstract}
\thanks{Research of the first author was partially supported by an NSERC
Discovery grant.}
\maketitle

\section{Introduction}
Let $S_k(N) $ denote the space of cusp forms of weight $k$  and level $N$ on the congruence subgroup ${\Gamma}_{0} (N)$ of $SL_2(\mathbb{Z})$. Let $f(z)= \sum_{n\geq 1} a_n q^n \in S_k(N)$ be a normalized Hecke eigenform, i.e., $f$ is an \emph{eigenfunction} for all the Hecke operators {$T_p$}'s and {$U_p$}'s and $a_1 = 1$.
The famous Ramanujan-Petersson conjecture (proved by Deligne \cite{deligne}) says that the $p$-th Fourier coefficient of $f$ is of the form
 
\begin{equation}\label{angle}
a_p (f)= 2 p^{(k-1)/2} \cos {\theta_p}(f), \quad p \not| N
\end{equation}
for some real \emph{angle} $\theta_p(f) \in [0, \pi]$.

\medskip

\noindent
In this paper we shall discuss the following question: 
\textit{for a fixed prime $p$, count the pairs of normalized
eigenforms $ (f,g)$ of weight $k$ and level $N$ such that $ a_p (f) = a_p (g).$}

\medskip

\noindent
More precisely, we prove the following:

\begin{theorem}\label{thm-1} 
For a fixed prime $p$, the number of pairs  $(f,g)$  
of normalized eigenforms in $S_k(N)$ such that $ \theta_p (f) = \theta_p (g)$ is bounded by $$ O \left( \frac{ \left( \dim S_k(N) \right)^2 (\log p )}{ (\log kN ) }\right),  $$
where the implied constant is absolute and independent of $p$.
\end{theorem}

In other words, our result gives a small saving over the trivial bound
provided $\log p \ll \log kN$.

Our interest in this question is partly motivated by a famous conjecture
of Maeda \cite{maeda} which predicts for $N=1$ that the polynomial
\begin{equation}\label{poly}
\prod_{f} (X - a_p(f)), 
\end{equation}
where the product is over all normalized Hecke eigenforms, is irreducible
over the rational number field.  In fact, Maeda conjectures that the Galois
group of this polynomial is the full symmetric group ${\frak S}_d$ where
$d$ is the dimension of the space $S_k(1)$.  In other words, Maeda
predicts that for level 1, the number of pairs in our theorem is
exactly dim $S_k(1)$.  For a fixed higher level, Tsaknias \cite{Tsak}
has conjectured that the above polynomial is a product of a bounded
number of irreducible polynomials viewed as a function of $k$.  
Though there is some computational data to support these conjectures,
they seem to be far out of reach of present knowledge and techniques.
Thus, it seems appropriate to investigate these questions through
methods currently known.  This is partial motivation of our work.

In the course of our proof, a certain exponential sum arises. Based on general heuristics about such exponential sums (more precisely 
the ``philosophy'' of \emph{square root} cancellation) we give a heuristic argument to support Maeda-Tsaknias's prediction \cite{Tsak}
regarding the bounded number of Galois orbits. 
Before we proceed to prove the theorem, we recall some preliminary results which will play an important role in proving the theorem.

\section{Approximation of characteristic functions with Selberg polynomials}  

Let $I = [a, b]$ be an interval contained in $[-1/2, 1/2]$
and $\chi_I$ the characteristic function
of the interval $I$. From the works of Selberg, Beurling
and Vaaler  (see \cite{vaaler}), there is a trigonometric polynomial $S_{M}$ of degree at most $M$ such that the following hold:
\begin{enumerate}
\item[(a)]
\begin{equation}\label{vaaler}
 \chi_I (x) \leq S_{M} (x),
\end{equation}
\item[(b)]
\begin{equation}\label{b}
\int_{-2}^{2} S_{M} (x) dx = b-a + \frac{1}{M+1}.
\end{equation}

\end{enumerate}

\noindent
Let $e(t)$ denote $e^{2\pi i t}$.
If we write the Fourier series for ${S}_{M} (x) $ as
\begin{equation}\label{fourier}
S_{M} (x) = \sum_{|n| \leq M} {\widehat{S}}_{M} (n) e(nx)
\end{equation}
then
\begin{enumerate}
\item[(c)]
\begin{equation}\label{basic}
|   {\widehat{S}}_{M} (n) | \leq \frac{1}{M+1} + \min \left( b-a , \frac{1}{\pi |n|} \right)
\end{equation}
\end{enumerate}

\section{A preliminary estimate}

From now onwards we shall denote by $\chi_I$ the characteristic function of the interval $ I = [ -\delta , \delta ] \subseteq  [ -1/2 , 1/2 ],$
with $\delta$ to be chosen later. 
\par
From \eqref{vaaler} it follows that
\begin{equation}\label{ineq-1}
\chi_I \left(\theta_p(f) - \theta_p(g) \right) \leq S_{M} \left(\theta_p(f) - \theta_p(g)\right).
\end{equation}
By the Fourier series expansion of  $ S_{M} (x)$, \eqref{ineq-1} can be written as
\begin{equation}\label{ineq-2}
\chi_I \left(\theta_p(f) - \theta_p(g) \right) \leq \sum_{-M}^{M} {\widehat{S}}_{M} (n) e\left(n \left(\theta_p(f) - \theta_p(g) \right) \right)
\end{equation}

\noindent
Thus the number of normalized Hecke eigenforms $f,g \in S_k(N)$ such that $\theta_p(f) = \theta_p(g)$  is
$$
\leq   \sum_{f\neq g} \chi_{I} ( \theta_p (f) - \theta_p (g) ) +\dim S_k (N)
=   \sum_{f, g} \chi_{I} ( \theta_p (f) - \theta_p (g) )
$$
For reasons that will become apparent later, it is convenient to 
write this as
$$\leq {1\over 2} \sum_{f,g} \chi_I(\pm \theta_p(f) \mp \theta_p(g)) $$
since $\chi_I$ is an even function. 
By \eqref{ineq-2}, this is 
$$
\leq {1 \over 2} \sum_{-M}^{M} {\widehat{S}}_{M} (n) \sum_{f,g} e\left(n \left(\pm \theta_p(f) \mp \theta_p(g) \right) \right) 
$$
This is  bounded by 
\begin{equation}\label{ineq-3}
 \leq {1 \over 2} 
\sum_{|n|\leq M} \mid {{\widehat{S}}_M} (n) \mid   {\left\vert \sum_f  
e ( \pm
n  \theta_p (f) ) \right \vert}^2 
\end{equation}
\noindent
Recall the estimate \eqref{basic}:
\begin{equation*}
| {\widehat{S}}_{M} (n) | \leq \frac{1}{M+1} + \min \left( b-a , \frac{1}{\pi |n|} \right)
\end{equation*}
Using this, the expression \eqref{ineq-3} for $n=0$ and for $n\neq 0$ (respectively) is
  \begin{equation}\label{key}
\begin{array}{lll}
&  \leq & \left( 2\delta +  \frac{1}{M+1}  \right)  {\left(  \dim S_k (N)  \right)}^2 \qquad  \\
 & & \\
& \qquad + &  \sum_{1 \leq |n| \leq M} \left\{ \frac{1}{M+1} + \min \left(2\delta, \frac{1}{\pi |n|} \right) \right\} {\left\vert \sum_f  e ( \pm  n\theta_p (f) )  \right\vert}^2 \\
\end{array}
\end{equation}

The crucial exponential sum that was elucidated in the beginning is $$  \sum_f  e ( \pm  n\theta_p (f) ) = \sum_{f} 2\cos n\theta_p(f)    .$$
It is this sum that appears in the Eichler-Selberg trace formula
and is thus amenable to estimation.  
In the next sections, we examine ways to estimate this sum.
\section{A heuristic argument}

We shall first give a heuristic argument to show that the estimate in the theorem is bounded by $O ( \left( {\dim S_k (N)} \right). $ 
In his 1997 paper \cite{Serre},  Serre  proved that the {$\theta_p$}'s are \emph{equidistributed} as $f$ varies with respect to a $p$-\emph{Sato-Tate measure}.
 This was made effective with error term by Ram Murty and Kaneenika Sinha in \cite{K-R}
using the Eichler-Selberg trace formula. 
More precisely, they gave precise estimates for
$$ \sum_{f} 2\cos n \theta_p(f) - c_n \dim S_k(N) , $$
where $c_0=1$, and $c_n= p^{-n/2} - p^{-(n-2)/2}$ for
$n$ even and zero  if $n$ is odd.    
It may be reasonable to expect \emph{square root} cancellation
for the error term, that is, 
$$
{\Big | \sum_f  e (\pm n  \theta_p (f) ) - c_n \dim S_k(N) \Big | } = O \left( \left( \dim S_k(N) \right)^{1/2} \right).
$$

Using the elementary inequality
$2ab \leq a^2 + b^2$, we have $|a+b|^2 \leq 2(|a|^2 + |b|^2)$, so that
$$\left| \sum_{f} e(\pm n\theta_p(f)) \right|^2 \ll |c_n|^2 ({\rm dim\, }S_k(N))^2 
+ \Big | \sum_f  e (\pm n  \theta_p (f) ) - c_n \dim S_k(N) \Big |^2. $$
With our assumption of square root cancellation, and taking $\delta = 1/M$, equation (\ref{key}) 
is
\begin{equation}
 \ll 
\frac{(\dim S_k(N))^2}{M+1} + 
\sum_{ \substack{1 \leq |n| \leq M}} \left\{ \frac{1}{M+1} + \min \left(2\delta, \frac{1}{\pi |n|} \right) \right\} \dim S_k (N) \end{equation}
\begin{equation*}
 \qquad  \le  \frac{(\dim S_k(N))^2}{M+1} + 
\frac{2M}{M+1} \dim S_k (N) ,
\end{equation*}
by virtue of the convergence of $$\sum_n c_n^2. $$  
Thus, by taking $M=\dim S_k(N)$, we obtain a final estimate of
\begin{equation}
  O \left( {\dim S_k (N)} \right).
\end{equation}

This estimate is consistent with the conjectures of Maeda and
Tsaknias.  Indeed, in the case of level 1, Maeda predicts that the
polynomial (\ref{poly}) is irreducible and so it should have {\bf no} repeated roots.
Thus
the number of pairs $(f,g)$ such that $\theta_p(f)=\theta_p(g)$
is equal to {\rm dim} $S_k(1)$.  In the higher level
case, Tsaknias predicts the number of Galois orbits is finite
and so the polynomial (\ref{poly}) is product of a {\bf finite}
number of irreducible polynomials with rational coefficients.
So in this case also, we expect the number to be bounded by
a constant multiple of {\rm dim} $S_k(N)$.

\section{Proof of the main theorem} 

What can be proved unconditionally? Recall that we want to find a non-trivial bound for the exponential sum $ \mid {\sum_f e ( \pm n \theta_p (f) ) \mid}^2 $ appearing earlier. This is achieved by using the following result of Ram Murty and Kaneenika Sinha (see \cite{K-R}, Theorem 18 ) which we state as

\begin{lemma} \label{lemma-1}
Define $f(N)$ as $\sum_{c|N} \phi({\rm gcd}(c,N/c)), $
and denote by $\nu(N)$ the number of distinct prime divisors of $N$.
Let
$$\psi(N) = N \prod_{p|N} \left( 1 + \frac{1}{p}\right). $$
Let $c_0=1$ and for $m \geq 1$, let
\[
 c_m = \begin{cases} p^{-m/2} - p^{-(m-2)/2} & \textrm{if\, } m\,  \textrm{is even}\\
0 & \textrm{if\, } m\,  \textrm{is odd.}\end{cases} \]
Then, 
\begin{equation}
| \sum_{f} 2 \cos m\theta_p(f) - c_m \text{\rm dim }S_k(N) |\leq 
4p^m 2^{\nu(N)} {\rm sup}_{f^2<4p^m} \psi(f) + 2f(N) + \delta_m(k), 
\end{equation}
where $\delta_m(k)=0$ unless $k=2$ in which case it is equal to
$2p^{m/2}$.

\end{lemma}
As mentioned in (11) of \cite{K-R}, the bound in the previous lemma can be
replaced by $p^{3m/2}2^{\nu(N)} \log p^m + \sqrt{N}d(N)$.
We have proved above that the quanity in question is given by (\ref{key}).
We choose $\delta = 1/M$.  Then our quantity is 
$$\ll {(\text{\rm dim }S_k(N))^2\over M} + \sup_{1\leq |m|\leq M}
\Big | \sum_f e(\pm m\theta_p(f))\Big |^2 $$
which is 
$$\ll {(\text{\rm dim }S_k(N))^2\over M} \left( 1+ \sum_{m\geq 1} {1\over p^m}
\right)
 + p^{3M}M^2 (\log p)^2 4^{\nu(N)} + Nd^2(N).$$
We need to make an optimal choice of $M$.  As the referee suggests,
this is best done by using the Lambert $W$-function.
We refer the reader to \cite{borwein} (see
especially Appendix A) for a friendly
introduction to this function.
Recall that this function is defined by
$$W(x)e^{W(x)} = x. $$
Our choice of $M$ will be so that the first two terms in the above
estimate are comparable.  That is, we seek $M$ so that
$$ ({\rm dim\, }S_k(N))^2 = p^{3M}M^3 (\log p)^2.$$
Since dim $S_k(N)$ lies somewhere between $kN$ and $kN \log \log N$,
our $M$ should (essentially) satisfy
$$ (kN)^2(\log p) = p^{3M} M^3 (\log p)^3. $$
Thus, 
$$M\log p  = W\left( (kN)^{2/3} (\log p)^{1/3}\right). $$
The Lambert function satisfies 
$$W(x) = \log x - \log \log x + o(1) $$
as $x$ tends to infinity, we deduce that $M$ can be taken as
the nearest integer to
$$  {{2\over 3} \log kN + {1 \over 3} \log p \over \log p} $$
which gives us a final estimate of
$$ \ll {(\textrm{dim }S_k(N))^2 \log p \over \log kN}. $$
This completes the proof.

\section{Concluding Remarks}

There are other predictions of the conjectures of Maeda and Tsaknias.
For instance, both conjectures predict that for a fixed $N$,
the number of normalized Hecke eigenforms with integers coefficients
is at most 1 (in the level one case) and $O_N(1)$ in the higher
level case.  In the latter case, the constant is expected to be
independent of $k$.  Such questions were investigated in \cite{K-R}
where similar methods and estimates were derived.   
$$\quad $$
\noindent {\bf Acknowledgements.}  We thank Sudhir Pujahari and
the referee for their
comments and corrections on an earlier draft of this paper.

\end{document}